\documentstyle[12pt]{article}
\begin{document}

\font\oo = cmr9
\font\nn = cmssbx10 
\font\fhead = cmbx10 at 12pt
\font\fhsl = cmbxsl10 at 12pt
\font\flp = cmmi12
\font\flpb = cmmib10 at 12pt
\font\flpp = cmmi10

\def\norm#1{$\vert\hskip-1.5pt\vert  \,#1\,
             \vert\hskip-1.5pt\vert$}
\def\nnorm#1{$\vert\hskip-1.5pt\vert\hskip-1.5pt\vert  \,#1\,
              \vert\hskip-1.5pt\vert\hskip-1.5pt\vert$}
\def\Norm#1{$\vert\hskip-2.9pt\vert\hskip-2.9pt\vert\hskip-2.9pt\vert \,#1\,
             \vert\hskip-2.9pt\vert\hskip-2.9pt\vert\hskip-2.9pt\vert$}

\def\hz{\smallskip}
\def\ez{\medskip}
\def\zz{\bigskip}
\def\dz{\medskip\bigskip}
\def\vz{\bigskip\bigskip}
\def\sz{\dz\dz}
\def\az{\vz\vz}

\def\noi{\noindent}
\def\LL{\Lambda}
\def\el{\,\raise 1pt\hbox{$\scriptstyle\in$}\,}
\def\mal{\raise 1.1pt\hbox{$\,\scriptscriptstyle\bullet\,$}}
\def\lp{\hbox{\flp `\hskip 0.8pt}}
\def\lpb{\hbox{\flpb `\hskip 0.8pt}}
\def\lpp{\hbox{\flpp `}}
\def\vre#1#2{\hskip #2pt\vrule width 0pt height #1 pt depth 0pt}
\def\vrf#1#2{\vrule width 0pt height #1 pt depth #2 pt}
\def\char#1{1_{\vrule width 0pt height 7.5pt depth 0pt #1}} 
\def\leer{\hskip2pt\not\hskip-2.5pt\raise 1.2pt\hbox{$\scriptstyle\bigcirc$}}
\def\Linf{$\hbox{\fhsl L}^\infty$}
\def\Lpinf{$\hbox{\fhsl L}^{p,\infty}$}
\def\Leins{$\hbox{\fhsl L}^1$}
\def\Lq{$\hbox{\fhsl L}^q$}
\def\index#1#2#3#4{\hbox{\flpp `}^#1_#2(\hbox{\flpp `}^#3_#4)}

\centerline{\large\bf THE $\hbox{K}_{\hbox{t}}$--FUNCTIONAL FOR THE}
\ez
\centerline{\large\bf INTERPOLATION COUPLE
$\hbox{L}^\infty$($d\mu\,;$\,$\hbox{L}^1$($d\nu$))\,,
$\hbox{L}^\infty$($d\nu\,;$\,$\hbox{L}^1$($d\mu$))}

\dz
\zz
\centerline{\large Albrecht He\ss}

\hz
\centerline{\large Gilles Pisier\footnote[1]
{\normalsize\vrf{13}{0}~Partially supported
by the N.S.F.}}

\az

\noi
{\bf \S 0. Introduction}
\dz

\noi
Let $(A_0,\,A_1)$ be a compatible couple of Banach spaces in the sense of
Interpolation Theory (cf. e.g. [BL]). Then the $K_t$--\,functional of an 
element $x$ in $A_0+A_1$ is defined for all $t>0$ as follows
$$K_t(x\,;\,A_0\,, A_1)=\inf\,
\bigl\{\Vert x_0\Vert_{A_0}+t\,\Vert x_1\Vert_{A_1}~~
\vert ~~x=x_0+x_1\bigr\}\,.$$

\noi
For instance the $K_t$--\,functional of the couple 
$(L^1({\bf R}), L^\infty({\bf R}))$ is well
known. We have
$$
K_t(x\,;\,L^1({\bf R}), L^\infty({\bf R}))=\int_0^t\,x^*(s)\,ds\\
=\sup\,\Bigl\{\int_E\,\vert x(s)\vert\,ds\,,~\vert E\vert=t\Bigr\}\,.
$$

\noi
See [BL] for more details and references. 

\hz
Let $(M,\mu)$ and $(N,\nu)$ be measure spaces. In this paper, we study the
$K_t$--\,functional for the couple
$$A_0=L^\infty(d\mu\,; L^1(d\nu))\,,~~A_1=L^\infty(d\nu\,; L^1(d\mu))\,.
\leqno(0.1)$$

\noi
Here, and in what follows 
the vector valued $L^p$--\,spaces $L^p(d\mu\,; L^q(d\nu))$ are meant 
in Bochner's sense.

\hz
One of our main results is the following, which can be viewed as a refinement
of a lemma due to Varopoulos [V].

\hz
\proclaim Theorem 0.1. Let $(A_0,A_1)$ be as in (0.1). Then for all $f$
in $A_0+A_1$ we have
$${1\over 2}\,K_t(f;\,A_0\,,A_1)\leq
\sup\,\bigg\{ \Big(\mu(E)\vee  t^{-1}\nu(F)\Big)^{-1}
\int_{E\times F} \vert f\vert\,d\mu\,d\nu\,\bigg\}
\leq K_t(f;\,A_0\,,A_1)\,,$$
where the supremum runs over all measurable subsets $E\subset M\,,~
F\subset N$ with positive and finite measure and $u\!\vee\!v$ denotes
the maximum of the reals $u$ and $v$.\vrf{0}{8}

\goodbreak

This result is a particular case of Theorem 3.2. below and its corollaries,
which give an analogous estimate for the couple
$$A_0=L^{p,\infty}(d\mu\,;\,L^q(d\nu))\,,~~
A_1=L^{p,\infty}(d\nu\,;\,L^q(d\mu))\,,$$

\noi
when $1\leq q<p\leq\infty\,.$ The preceding Theorem 0.1 corresponds to the
case $p=\infty,~q=1$.

\hz
The paper is organized as follows. In \S 1, we prove Theorem 0.1 
in the finite
discrete case, i.e. when $M$ an $N$ are finite sets equipped with discrete
measures. In \S 2, we discuss the extension to 
$L^{p,\infty}(d\mu\,;\,L^q(d\nu))\,,~L^{p,\infty}(d\nu\,;\,L^q(d\mu))\,,$
again in the discrete case and treat the case of general measure spaces
in \S 3. Finally in \S 4 we give some applications.

\hz
To describe the main one, let us denote by $B_1$ (resp. $B_0$) the space
of all bounded operators from $L^1(d\nu)$ to $L^1(d\mu)$ 
(resp. from $L^\infty(d\nu)$ to $L^\infty(d\mu)$). Our results yield 
a description of the space $(B_0,B_1)_{\theta,q},
~0<\theta<1,~1\leq q\leq\infty$, obtained by the real (i.e. Lions\,--\,Peetre)
interpolation method. The result is particularly simple  in the case
$q=\infty$. In that case, we prove

\hz
\proclaim Theorem 0.2. The space $(B_0,B_1)_{\theta,\infty}$ coincides
with the space of all bounded regular operators $u$ from $L^{p,1}(d\nu)$
to $L^{p,\infty}(d\mu)$ with $p=1/\theta$.

\noi
Here ''regular`` is meant in the sense of [MN]. Equivalently $u$ belongs
to $(B_0,B_1)_{\theta,\infty}$ iff there is a positive (i.e. positivity
preserving) operator $v$ from $L^{p,1}(d\nu)$
to $L^{p,\infty}(d\mu)$ which dominates $u$, i.e. such that $-v\leq u\leq v.$
(Note that we implicitly use only real scalars, but this is not essential.)
In particular, it follows that $u$ belongs to 
$(B_0,B_1)_{\theta,\infty}$ iff $\vert u\vert$ is bounded from
$L^{p,1}(d\nu)$ to $L^{p,\infty}(d\mu)$, or equivalently iff $\vert u\vert$
is of ''very weak type $(p\,,p)$`` in the sense of [SW, chapter 3.3].
We refer to [MN] for the definition of the modulus $\vert u\vert$ of
a regular operator acting between Banach lattices. We merely recall
that if $u$ is given by a matrix $(a_{ij})$ acting between two
sequence spaces, then $\vert u\vert$ corresponds to the matrix
$(\vert a_{ij}\vert).$ A similar fact holds with ''kernels`` instead
of matrices.

\vfill
\eject
\noi
{\bf \S 1. The $\hbox{\bf K}_{\hbox{\bf t}}$--\,functional for the 
interpolation couple            $\hbox{\lpb}_m^\infty(\hbox{\lpb}_n^1)\,,
                   \hbox{\lpb}_n^\infty(\hbox{\lpb}_m^1)\,$}
\dz

\noi We use the following notation:
for $p\el{\bf R}\,,\,p\!\geq\!1\,,$ let $p^\ast$ 
be the conjugate exponent of $p\,,$ $1/p\!+\!1/{p^\ast}\!=\!1$, the dot 
\mal\ denotes the pointwise multiplication of matrices $a\,,b$ according to
$(a\mal b)(i,j)\,=\,a(i,j)b(i,j)\,.$  For a 
set $A$ let $\char{A}$ be the characteristic function of $A$ 
(the whole set will always be known from the context).

\hz

Let $(M,\mu)$ be the measure space consisting of the atoms $\{1,\dots,m\}$ 
of positive masses $\mu_1,\dots,\mu_m$ and $(N,\nu)$ the measure 
space consisting of the atoms $\{1,\dots,n\}$ of positive masses 
$\nu_1,\dots,\nu_n$. We equip $M\!\times\!N$ with the product measure
$\mu\!\times\!\nu$.

\hz

Besides the $\hbox{\lp}^p\!$--norms 
$$\hbox{\norm{a}}_{\lpp^p} = \Bigl(\;\sum_{(i,j)\in M\times N}
                      \vert a(i,j)\vert^p\mu_i\nu_j\Bigr)^{1/p}\,,$$

\noi we introduce in this section for $m\!\times\!n$--matrices
$a$ the norms

$$\hbox{\norm{a}}
     =   \hbox{\norm{a}}_{\hbox{\flpp `}^\infty_m(\hbox{\flpp `}^1_n)}
     =   \max_{i\in M}\;\Bigl(\sum_{j\in N}\vert a(i,j)\vert\,\nu_j\Bigr)\,,
$$

$$\hbox{\norm{a}}{\vre{3}{-0.1}}^{\top} 
     =   \hbox{\norm{a^\top}}_{\hbox{\flpp `}^\infty_n(\hbox{\flpp `}^1_m)}
     =   \max_{j\in N}\;\Bigl(\sum_{i\in M}\vert a(i,j)\vert\,\mu_i\Bigr)\,.
$$       
                 
\noi It is straightforward to see that for $E\subset M\,,\,F\subset N\,,$
$$\hbox{\norm{\char{E\times F}\mal a}}_{\lpp^1} 
   \leq \min\; \bigl\{ \mu(E)\hbox{\norm{a}},\, 
                     \nu(F){\hbox{\norm{a}}}{\vre{3}{-0.1}}^{\top}
             \bigr\}\,.
$$
\hz
\noi Therefore we have for matrices $b,\,c$ and any $t>0$
$$\hbox{\norm{\char{E\times F}\mal (b+c)}}_{\lpp^1} 
 \leq \Bigl(\mu(E)\lor t^{-1}\nu(F)\Bigr)
      \Bigl(\hbox{\norm{b}}+t\,\hbox{\norm{c}}{\vre{3}{-0.1}}^{\top}\Bigr)\,.
\leqno(1.1)$$

\ez

If we introduce the norm
$$\hbox{\nnorm{a}}_t = \sup_{E,F} \,
\Bigl(\mu (E)\lor t^{-1}\nu (F)\Bigr)^{-1}
\hbox{\norm{\char{E\times F}\mal a}}_{\lpp^1}\,,  
$$

\noi then  (1.1) means that the $K_t$--\,functional
for the interpolation couple 
$\hbox{\lp}_m^\infty(\hbox{\lp}_n^1)\,,
\hbox{\lp}_n^\infty(\hbox{\lp}_m^1)\,$  
$$K_t(a) = \inf\bigl\{\hbox{\norm{b}}+t\,\hbox{\norm{c}}
{\vre{3}{-0.1}}^{\top}\,\vert~\,
a=b+c\,\bigr\}$$

\noi is an upper bound for this norm
$$\hbox{\nnorm{a}}_t \leq K_t(a)\,.$$
\hz

\noi In the following we prove an upper estimate of $K_t(\,.\,)$ by
the functional $\hbox{\nnorm{.}}_t\,.$

\hz

\proclaim Proposition 1.1. Let $t>0$. For any matrix $a$ with 
$\hbox{\nnorm{a}}_t\leq 1$ we can
split $M\!\times\!N$ into a disjoint union $M\!\times\!N = A\cup B\,,\ 
A\cap B=\leer\,,$ such that
$$\hbox{\norm{\char{A}\mal a}} \leq 1\,,~~~
  \hbox{\norm{\char{B}\mal a}}{\vre{3}{-0.1}}^{\top} \leq 1/t\,,
\leqno(1.2)$$
hence
$$K_t(a) \leq 2\,\hbox{\nnorm{a}}_t \,.$$
\par
\medskip
\noi{\bf Proof.} We proceed by induction on $m$. For $m\!=\!1$ 
suppose, without loss of generality, that
$$\vert a(1,1)\vert\geq\vert a(1,2)\vert\geq\dots\geq\vert a(1,n)\vert\,.
\leqno(1.3)$$

\noi If $t^{-1}\nu(N)\!\leq\!\mu(\{1\})\!=\!\mu_1$ 
we put $k=n$, if not, let $k\,,~0\!\leq\!k\!<\!n\,,$ 
be as large as possible 
such that $t^{-1}\nu( \{1,\dots,k\})\leq\mu_1$. For 
$E=\{1\}$ and $F=\{1,\dots,k\}$ we have in any case
$$\Bigl(\mu (E)\lor t^{-1}\nu (F)\Bigr)^{-1}
  \hbox{\norm{\char{E\times F}\mal a}}_{\lpp^1} =
\sum_{j\leq k}\vert a(1,j)\vert\,\nu_j\leq 1.
\leqno(1.4)$$

\noi If $k<n$ then $t^{-1}\nu(\{1,\dots,k\!+\!1\})>\mu_1$, so we
obtain for 
$E=\{1\}$ and $F=\{1,\dots,k\!+\!1\}$ 

$$\Bigl(\mu (E)\lor t^{-1}\nu (F)\Bigr)^{-1}
  \hbox{\norm{\char{E\times F}\mal a}}_{\lpp^1} = \mu_1 t~
{{\sum_{j\leq k+1}{\vrule width 0pt height 0pt depth 5pt}
\vert a(1,j)\vert\,\nu_j}
\over{\sum_{j\leq k+1}\nu_j}}  
\leq 1.$$

\noi By (1.3) this yields 
$$\vert a(1,n)\vert\,\mu_1\leq\,\dots\,\leq\vert a(1,k\!+\!1)\vert\,\mu_1\leq
1/t\,.\leqno(1.5)$$

\noi Take now
$A=\{1\}\!\times\!\{1,\dots,k\}\,,~B=(M\!\times\!N)\!\setminus\!A\,,$
\noi and (1.2) is fulfilled.
\ez

\noi Let us assume the truth of the proposition for $1,\dots,m\!-\!1$. 
Without loss of generality, we can suppose the sums over $M$
$$\sigma_j=\sum_{i\leq m}\vert a(i,j)\vert\,\mu_i$$

\noi to be in descending order 
$\sigma_1\geq\sigma_2\geq\dots\geq\sigma_n$. If $t^{-1}\nu(N)\leq\mu(M)$ let
$k=n$, if not, let $k\,,~0\!\leq\!k\!<\!n\,,$ be the largest number such that
$t^{-1}\nu(\{1,\dots,k\})\leq\mu(M)$. For 
$E\!=\!M$ and $F\!=\!\{1,\dots,k\}$ we have
$$\hbox{\norm{\char{E\times F}\mal a}}_{\lpp^1}\leq\mu(M).
\leqno(1.6)$$

\noi If we permute the rows such that  
$$\tau_i=\sum_{j\leq k}\vert a(i,j)\vert\,\nu_j$$

\noi are in descending order $\tau_1\geq\tau_2\geq\dots\geq\tau_m$ then
it follows from (1.6) that
$$\tau_m\leq 1.\leqno(1.7)$$

\noi By the induction hypothesis on $\{1,\dots,m\!-\!1\}\!\times\!
\{1,\dots,k\}$ we find a splitting
$$\{1,\dots,m\!-\!1\}\!\times\!\{1,\dots,k\}=A_1\cup B_1\,,~~
A_1\cap B_1=\leer\,,$$

\noi such that
$$\hbox{\norm{\char{A_1}\mal a}} \leq 1\,,~~~
  \hbox{\norm{\char{B_1}\mal a}}{\vre{3}{-0.1}}^{\top} \leq 1/t\,.
$$

\noi For the wanted splitting of $M\!\times\!N$
we only have to put
$$A=A_1\cup (\{m\}\!\times\!\{1,\dots,k\})\,,~~~
B=(M\!\times\!N)\!\setminus\!A.$$

\noi Indeed, by the induction hypothesis on $A_1$ and (1.7)
$$\hbox{\norm{\char{A}\mal a}} \leq 1\,,$$

\noi and we are done for $k\!=\!n$ which implies $B\!=\!B_1\,.$
If $k<n$ we argue as in
 case $m=1$ (using $E=M$ and $F=\{1,\dots,k\!+\!1\}$) and show that 
$$\sigma_n\leq\dots\leq\sigma_{k+1}\leq 1/t\,.$$

\noi This yields together with the hypothesis on $B_1$
$$\hbox{\norm{\char{B}\mal a}}{\vre{3}{-0.1}}^{\top}  \leq 1/t\,.$$

\hz

\noi {\bf Remark 1.2.} There is another ``norm''
$$\hbox{\Norm{a}}_t = \sup\,\bigl\{
\mu(E)^{-1} \hbox{\norm{\char{E\times F}\mal a}}_{\lpp^1}\,~
\vert\,~t^{-1}\nu(F)\leq\mu(E)\bigr\}  
$$

\noi closely related to $\hbox{\nnorm{a}}_t\,.$ This functional
$\hbox{\Norm{a}}_t$ has the drawback that for small values of $t$ or for
strongly varying masses the supremum may be on an empty set, 
in that case we put $\hbox{\Norm{a}}_t =0\,.$
In any case,
$$\hbox{\Norm{a}}_t\leq\hbox{\nnorm{a}}_t\leq K_t(a)\,.$$

\noi  But for uniform masses, say $\mu_i\!=\!\nu_j\!=\!1\,,$ and for 
$t\!\geq\!1\,,$ where $\hbox{\Norm{.}}_t$ is indeed a norm,
one can construct for a matrix $a$ with $\hbox{\Norm{a}}_t\leq 1$ 
almost along the same lines as above a splitting 
$M\!\times\!N = A\cup B\,,\ A\cap B=\leer\,,$ such that
(denoting by $[\,t\,]$ the integer part of $t$)
$$\hbox{\norm{\char{A}\mal a}} \leq 1\,,~~~
\hbox{\norm{\char{B}\mal a}}{\vre{3}{-0.1}}^{\top} \leq 1/[\,t\,]\>,$$
\noi hence 
$$K_t(a)\leq (1+t/[\,t\,])\hbox{\Norm{a}}_t<3\,\hbox{\Norm{a}}_t\,,$$
\ez
\noi an inequality due, for $t\!=\!1,$ to Varopoulos [V], cf. also [BF].

\zz

\noi
{\bf \S 2. The $\hbox{\bf K}_{\hbox{\bf t}}$--\,functional for the 
interpolation couple       
       $\hbox{\lpb}_M^{p,\infty}\!(\hbox{\lpb}_N^q)\,,
             \hbox{\lpb}_N^{p,\infty}\!(\hbox{\lpb}_M^q)$     }
\dz

\noi For any Bochner--measurable function $f$ on a measure space 
$(\Omega,\Sigma,\mu)$ with values in a Banach space $X$ and
for $p>0$ let
$$\hbox{\norm{f}}_{L^{p,\infty}(X)} = 
\inf\,\bigl\{C\,\vert~t^p\mu\{\Vert f\Vert\!>\!t\}\leq C^p~~
\hbox{for all}~~t\!>\!0\bigr\}\,,
\leqno(2.1)$$

\noi or more generally by using the  nonincreasing, equimeasurable 
rearrangement $f^\ast$ of $\Vert f\Vert$
$$\hbox{\norm{f}}_{L^{p,\infty}(X)} =
\sup_{t>0}\,t^{1/p}f^\ast\!(t)\,,~~
\hbox{\norm{f}}_{L^{p,q}(X)} =
\Bigl(\int_0^\infty(t^{1/p}f^\ast(t))^q\,{dt\over t}\Bigr)^{1/q}\>.
\leqno(2.2)$$

\noi In general, these quantities are not norms but can be 
replaced for $p\!>\!1,\,q\!\geq\!1$ by equivalent norms. 
We do not use this and refer to [BS], Lemma IV.4.5, p.\,219. 

\hz

We start with the case of the interpolation couple       
$\hbox{\lp}_m^{p,\infty}(\hbox{\lp}_n^1)\,,
\hbox{\lp}_n^{p,\infty}(\hbox{\lp}_m^1)$ on the finite 
measure spaces $(M,\mu)$, $(N,\nu)$ from section 1
and consider the functionals
$$\hbox{\norm{a}}
     =   \hbox{\norm{a}}_{\hbox{\flpp `}^{p,\infty}_m(\hbox{\flpp `}^1_n)}~~
\hbox{and}
~~\hbox{\norm{a}}{\vre{3}{-0.1}}^{\top} 
  = \hbox{\norm{a^\top}}_{\hbox{\flpp `}^{p,\infty}_n(\hbox{\flpp `}^1_m)}\>.
$$

\noi We have for $p>1$ and $E\subset M\,,~F\subset N$
$$\hbox{\norm{\char{E\times F}\mal a}}_{\lpp^1} 
   \leq p^\ast\min\; \bigl\{ \mu(E)^{1/p^\ast}\!\hbox{\norm{a}}\,,~ 
                   \nu(F)^{1/p^\ast}\!{\hbox{\norm{a}}}{\vre{3}{-0.1}}^{\top}
             \bigr\}\,.
\leqno(2.3)$$

\hz

We introduce the  norm $\hbox{\nnorm{.}}_{p,t}$ as follows
$$\hbox{\nnorm{a}}_{p,t} = \sup_{E,F} \,
\Bigl(\mu (E)^{1/p^\ast}\!\lor t^{-1}\nu (F)^{1/p^\ast}\Bigr)^{-1}
\hbox{\norm{\char{E\times F}\mal a}}_{\lpp^1}\,.  
$$

\noi As in section 1 we want to compare the $K_t$--\,functional 
$K_{p,t}(\,.\,)$ for the interpolation couple 
$\hbox{\lp}_m^{p,\infty}(\hbox{\lp}_n^1)\,,
\hbox{\lp}_n^{p,\infty}(\hbox{\lp}_m^1)$ 
with the norm 
$\hbox{\nnorm{.}}_{p,t}\,.$ 
In view of (2.3) we have a lower estimate of $K_{p,t}(a)$ 
$$\hbox{\nnorm{a}}_{p,t}\leq p^\ast K_{p,t}(a)\,.
\leqno(2.4)$$

\noi Now we establish the corresponding upper estimate.

\hz

\proclaim Proposition 2.1. Let $t>0\,,~p>1$. For any matrix $a$ with 
$\hbox{\nnorm{a}}_{p,t}\leq 1$ we can
split $M\!\times\!N$ into a disjoint union $M\!\times\!N = A\cup B\,,\ 
A\cap B=\leer\,,$ such that
$$\hbox{\norm{\char{A}\mal a}} \leq 1\,,~~~
  \hbox{\norm{\char{B}\mal a}}{\vre{3}{-0.1}}^{\top} \leq 1/t\,,
\leqno(2.5)$$
hence
$$K_{p,t}(a) \leq 2\,\hbox{\nnorm{a}}_{p,t} \,.$$
\par
\medskip
\noi{\bf Proof.} Similar to the proof of Proposition 1.1.
For $m\!=\!1$ take $k\!=\!n$ if $\nu(N)\leq\mu_1\,t^{p^\ast}\,,$ if not,
suppose 
$\vert a(1,1)\vert\geq\vert a(1,2)\vert\geq\dots\geq\vert a(1,n)\vert$ 
and choose $k\,,~0\!\leq\!k\!<\!n\,,$ as large as possible
such that $\nu( \{1,\dots,k\})\leq\mu_1\,t^{p^\ast}\,.$ 
\vre{10}{0} Inserting $E=\{1\}$ and $F=\{1,\dots,k\}$  we have 
$$\mu(\{1\})^{1/p}\sum_{j\leq k}\vert a(1,j)\vert\,\nu_j\leq 1\,,$$

\noi and for 
$E=\{1\}$ and $F=\{1,\dots,l\}\,,~l>k\,,$ 
$$\nu(\{1,\dots,l\})^{1/p}\vert a(1,l)\vert\mu_1\leq 1/t\,.$$

\noi So we can put
$A=\{1\}\!\times\!\{1,\dots,k\}\,,~B=(M\!\times\!N)\!\setminus\!A$
\noi to obtain (2.5).
\ez
 
\noi For $m\!>\!1$ 
let us put $k\!=\!n$ if $\nu(N)\!\leq\!\mu(M)\,t^{p^\ast}\!,$ if not, 
choose the largest
$k\,,~0\!\leq\!k\!<\!n\,,$ with
$\nu( \{1,\dots,k\})\leq\mu(M)\,t^{p^\ast}$. Taking 
$E=M$ and $F=\{1,\dots,k\}$ we have
$$\hbox{\norm{\char{E\times F}\mal a}}_{\lpp^1}\leq\mu(M)^{1/p^\ast}\,,$$

\noi consequently
(we use the notations $\tau_m$, $\sigma_l\,$ of the preceding proof)
$$\mu(\{1,\dots,m\})^{1/p}\tau_m\leq 1\,,\leqno(2.6)$$

\noi and for $E=M$ and $F=\{1,\dots,l\}\,,~l>k\,,$
$$\nu(\{1,\dots,l\})^{1/p}\,\sigma_l\leq 1/t\,.\leqno(2.7)$$

\noi The induction hypothesis on
$\{1,\dots,m\!-\!1\}\!\times\!\{1,\dots,k\}$ together with (2.6) and (2.7)
yields the result.

\zz

\noi {\bf Remark 2.2.} The inequality (2.4) is meaningless for 
$p=1\,,$ but the corresponding norm 
$$\hbox{\nnorm{a}}_{1,t} = \sup_{E,F} \;
(1\!\land\!t)\;\hbox{\norm{\char{E\times F}\mal a}}_{\lpp^1}
= (1\!\land\!t)\; \hbox{\norm{a}}_{\lpp^1(M\times N)}$$

\noi is proportional to the norm 
$\hbox{\norm{a}}_{\lpp^1(M\times N)}\,.$

\zz

\noi {\bf Remark 2.3.} For 
$M\!=\!N\!=\!\{1,\dots,m\}\,,\,\mu_i\!=\!\nu_j\!=\!1\,,$ 
and for $t\geq 1$ one can prove with only slight changes the 
equivalence of $\hbox{\nnorm{.}}_{p,t}$ and
$$\hbox{\Norm{a}}_{p,t} = \sup\,\bigl\{
\mu(E)^{-1/{p^\ast}} \hbox{\norm{\char{E\times F}\mal a}}_{\lpp^1}\,~
\vert\,~t^{-1}\nu(F)^{1/p^\ast}\leq\mu(E)^{1/p^\ast}\bigr\}\,,$$

\noi namely, with constants independent of $m$ and $t$
$$(1/p^\ast)\,\hbox{\Norm{a}}_{p,t}\leq
(1/p^\ast)\,\hbox{\nnorm{a}}_{p,t}\leq
K_{p,t}(a)\leq 3\,\hbox{\Norm{a}}_{p,t}\,.$$

\noi We leave the details to the reader.
>From the example 
$M\!=\!\{1\}\,,$ $~N\!=\!\{1,\dots,n\}\,,$ the
$\nu_j$ and  $\mu_1$ equal to $1$ and $a\equiv 1\,,$ 
where
${\hbox{\nnorm{a}}}_{p,1} = 1\,,$ 
but $K_{p,1}(a)=n^{1/p}\,,$
we see that there should be a restriction $M\!=\!N\,$
in order to have the occurring  constants independent of $n$. 

\zz

\noi {\bf Remark 2.4.} If we replace $\hbox{\norm{a}}$ by
$\hbox{\norm{a}}_{\hbox{\flpp `}^{p,q}_m(\hbox{\flpp `}^1_n)}$
and $\hbox{\norm{a}}{\vre{3}{-0.1}}^{\top}$ by
$\hbox{\norm{a^\top}}_{\hbox{\flpp `}^{p,q}_n(\hbox{\flpp `}^1_m)}$
then
$$\hbox{\norm{\char{E\times F}\mal a}}_{\lpp^1} 
\leq C(p,q)\,\min\; \bigl\{ \mu(E)^{1/p^\ast}\!\hbox{\norm{a}}\,,~ 
\nu(F)^{1/p^\ast}\!{\hbox{\norm{a}}}{\vre{3}{-0.1}}^{\top}\bigr\}$$

\noi for the constant $C(p,q)\!=\!(p^\ast\!/q^\ast)^{1/q^\ast}.$
Therefore the $K_t$--func\-tio\-nal for the  
interpolation couple
$\hbox{\lp}_m^{p,q}(\hbox{\lp}_n^1)\,,\hbox{\lp}_n^{p,q}(\hbox{\lp}_m^1)\,$
is greater than $C(p,q)^{-1}\hbox{\nnorm{.}}_{p,t}\,.$ But there is
no upper bound of this {$K_t$--\,functional} 
\vre{10}{0} by $\hbox{\nnorm{.}}_{p,t}$ 
independent of $m$ and $n$.

\hz

\noi Take e.g. $M=N=\{1,\dots,n\}\,,~\mu_i=\nu_j=1$ and
$$a=\left(\matrix{1 & 2^{-1/p} & \ldots & n^{-1/p} \cr
                  0 &   0      & \ldots &   0      \cr
             \vdots &  \vdots  &        & \vdots   \cr          
                  0 &   0      & \ldots &   0      \cr}\right)$$

\noi then 
$\hbox{\nnorm{a}}_{p,t}\leq p^\ast\,,~K_t(a)\sim (\log n)^{1/q}\,.$

\ez\hz

If we now put for $1\leq p,q\leq\infty$  
$$\hbox{\norm{a}} =
\hbox{\norm{a}}_{\hbox{\flpp `}^{p,\infty}_m(\hbox{\flpp `}^q_n)}\,,~~~
\hbox{\norm{a}}{\vre{3}{-0.1}}^{\top} =
\hbox{\norm{a^\top}}_{\hbox{\flpp `}^{p,\infty}_n(\hbox{\flpp `}^q_m)}\,,
\leqno(2.8)$$

\noi define the $K_t$--\,functional $K_{p,q,t}(\,.\,)$ for 
the interpolation couple
$\hbox{\lp}_m^{p,\infty}(\hbox{\lp}_n^q)\,,
\hbox{\lp}_n^{p,\infty}(\hbox{\lp}_m^q)$
and let, with the abbreviation $\alpha\!=\!1/q\!-\!1/p\,,$
$$\hbox{\nnorm{a}}_{p,q,t} = \sup_{E,F} \,
\Bigl(\mu (E)^{\alpha}\!\lor t^{-1}\nu (F)^{\alpha}\Bigr)^{-1}
\hbox{\norm{\char{E\times F}\mal a}}_{\lpp^q}\,,
\leqno(2.9)$$

\noi then we derive from Proposition 2.1 by $q$-convexification 
(cf.\,\,[LT]) 
the following theorem.

\zz

\proclaim Theorem 2.5. Let $t>0\,,~1\leq q<p\leq\infty\,.$ Denoting
$C(p,q)\!=\!(1\!-\!q/p)^{1/q}$ we can estimate
the $K_t$--\,functional $K_{p,q,t}(\,.\,)$ for 
\vre{10}{0} the interpolation couple
$\hbox{\lp}_m^{p,\infty}(\hbox{\lp}_n^q)\,,
\hbox{\lp}_n^{p,\infty}(\hbox{\lp}_m^q)$
by
$$C(p,q)\,\hbox{\nnorm{a}}_{p,q,t}\leq K_{p,q,t}(a)
\leq 2\,\hbox{\nnorm{a}}_{p,q,t}\,.\leqno(2.10)$$\hz
More precisely, for $\hbox{\nnorm{a}}_{p,q,t}\leq 1$ we find $A\,,\,B$ with
$M\!\times\!N =A\cup B\,, ~A\cap B =\leer\,,$ and
$$\hbox{\norm{\char{A}\mal a}} \leq 1\,,~~~
  \hbox{\norm{\char{B}\mal a}}{\vre{3}{-0.1}}^{\top} \leq 1/t\,.$$
\par 
\ez 

\noi {\bf Remark 2.6.} An analogous result can be stated for 
measure spaces $M\!=\!N$
with equal atoms $\mu_i=\nu_j=1\,$ and $t\geq 1\,$ for the norm
$$\hbox{\Norm{a}}_{p,q,t} = \sup\,\bigl\{
\mu(E)^{-\alpha} \hbox{\norm{\char{E\times F}\mal a}}_{\lpp^1}\,~
\vert\,~t^{-1}\nu(F)^\alpha\!\leq\mu(E)^\alpha\bigr\}\,,~~~~
(\alpha\!=\!1/q\!-\!1/p) $$
\noi (cf. Remark 2.3).

\zz

\proclaim Theorem 2.7. The statement (2.10) remains true for arbitrary
discrete measure spaces. 
\par
\zz 
\noi{\bf Proof.} Let $M$, $N$ be arbitrary discrete measure spaces.
It is very easy to check that for an element $a=(a(i,j))$ we have
$$\Vert a\Vert_{\hbox{\lp}_M^{p,\infty}\!(\hbox{\lp}_N^q)}
= \sup\;\Vert 1_{A\times B}\mal a\Vert_
{\hbox{\lp}_M^{p,\infty}\!(\hbox{\lp}_N^q)}\,,$$

\noi where the sup runs over all finite subsets $A\subset M$,
$B\subset N$. Hence the result follows by a simple pointwise
compactness argument left to the reader.

\vfill
\eject

\noi
{\bf \S 3. The $\hbox{\bf K}_{\hbox{\bf t}}$--\,functional for the
couple 
\Lpinf ($d\mu\,;$\,\Lq ($d\nu$))\,, \Lpinf ($d\nu\,;$\,\Lq ($d\mu$))}

\dz

\noi Now we generalize the results from the previous sections to 
arbitrary measure spaces and treat the generic case
$(M,\mu)\!=\!(N,\nu)\!=\!({\bf R},\lambda)$ 
of non--atomic measure spaces, where 
$\lambda$ is the Lebesgue\,--measure. 
For
$a\!=\!a(x,y)\el L_{loc}^q({\bf R}^2)$ let us define the functionals
$$\hbox{\norm{a}}
=   \hbox{\norm{a}}_{L^{p,\infty}_x(L^q_y)}\,,~~~~
\hbox{\norm{a}}{\vre{3}{-0.1}}^{\top} 
=   \hbox{\norm{a^\top}}_{L^{p,\infty}_y(L^q_x)}\>.$$

\noi These  
functionals are used to define the 
$K_t$--\,functional $K_{p,q,t}(\,.\,)$ of the interpolation couple
$L^{p,\infty}(dx\,;L^q(dy))\,,L^{p,\infty}(dy\,;L^q(dx))$
in the obvious manner.
For $t\!>\!0,$ 
$1\!\leq\!q\!<\!p\!\leq\!\infty$ and $\alpha\!=\!1/q\!-\!1/p,$ 
let us introduce the norm
\hz
$$\hbox{\nnorm{a}}_{p,q,t} = \sup_{E,F} \,
\Bigl(\lambda(E)^{\alpha}\!\lor t^{-1}\lambda(F)^{\alpha}\Bigr)^{-1}
\Bigl(\int\!\!\!\int_{E\!\times\!F}\vert a(x,y)\vert^q\,dx\,dy\Bigr)
^{1/q}\,.$$

\noi Obviously (c.f. Theorem 2.5) for $C(p,q)=(1-p/q)^{1/q}$
$$C(p,q)\,\hbox{\nnorm{a}}_{p,q,t}\leq K_{p,q,t}(a)\,.$$

\hz

\noi We will now prove the counterpart
$$K_{p,q,t}(a)\leq 2\,\hbox{\nnorm{a}}_{p,q,t}\,.$$

\proclaim Theorem 3.1. For any $t>0$ we have 
the inequality
$$C(p,q)\,\hbox{\nnorm{a}}_{p,q,t}\leq K_{p,q,t}(a)\leq
2\,\hbox{\nnorm{a}}_{p,q,t}
\leqno(3.1)$$
between the norm
$\hbox{\nnorm{.}}_{p,q,t}$ defined above
and the $K_t$--\,functional of the interpolation couple
$L^{p,\infty}(dx\,;L^q(dy))\,,L^{p,\infty}(dy\,;L^q(dx))\,.$ 
\par
 
\zz

\noi {\bf Proof.} Let us call countably simple a function
$a:{\bf R}^2\!\rightarrow {\bf R}$ of the form
$a=\sum a_{ij}\,1_{A_i\times B_j}$ where $(A_i)$ and $(B_j)$ are 
countable measurable partitions of {\bf R}.
One can easily check that the subset of countably
simple functions is dense in 
$L^{p,\infty}(dx\,;L^q(dy))+L^{p,\infty}(dy\,;L^q(dx))\,,$
hence it suffices to check the inequality (3.1) for those functions.
But in that case (3.1) follows from Theorem 2.7.

\noi (This argument presupposes that the functional 
$K_{p,q,t}$ is continuous with respect to the norm 
$\hbox{\nnorm{.}}_{p,q,t}$, but this can be checked using the
equivalent renorming of weak $L^p$, conditional expectations and
a weak compactness argument.)

\ez

\proclaim Corollary 3.2. By combination of Theorems 2.7 and 3.1 we
obtain the inequality (3.1) for all measure spaces.

\goodbreak
\ez

\noi {\bf Remark 3.3.} In the case of non--atomic measure
spaces $(M,\mu),$ and $(N,\nu)$ of infinite total measures
$\mu(M)\!=\!\nu(N)\!=\!\infty$ the difference between 
$\hbox{\nnorm{a}}_{p,q,t}$ and 
$$\hbox{\Norm{a}}_{p,q,t} = \sup\,\bigl\{
\mu(E)^{-\alpha} \hbox{\norm{\char{E\times F}\,a}}_{L^q(M\times N)}\,~
\vert\,~t^{-1}\nu(F)^\alpha\leq \mu(E)^\alpha\bigr\}~~~~
(\alpha\!=\!1/q\!-\!1/p)$$

\noi disappears completely. Obviously
$$\hbox{\nnorm{a}}_{p,q,t} = \hbox{\Norm{a}}_{p,q,t}\,.$$

\noi If the (non--atomic) measure spaces
$(M,\mu)\,,\,(N,\nu)$ have the same finite total measure, say $1$,
the equivalence of these two norms with constants independent
of $t$ can be established for $t\geq 1$. 
The example of the function
$a\equiv 1$ on $[\,0,1\,]\!\times\![\,0,1\,]$ (equipped with the
Lebesgue measure) gives for $0\!<\!t\!\leq\!1$
$$\hbox{\nnorm{a}}_{p,q,t}=t~~\hbox{and}~~
\hbox{\Norm{a}}_{p,q,t}=t^{p/p-q}\,.$$

\noi This leads for $0\!<\!t\!\leq\!1$ to the inequality
$$\sup_{a\in L^1(M\times N)}\,
{{\hbox{\nnorm{a}}_{p,q,t}}\over{\hbox{\Norm{a}}_{p,q,t}}}\geq
t^{-q/p-q}\,.$$

\noi But we have indeed

\zz

\proclaim Proposition 3.4.  For non--atomic measure spaces 
$(M,\mu)\,,\,(N,\nu)$ with the same finite total measure, say 
$\mu(M)\!=\!\nu(N)\!=\!1$,
$1\!\leq\!q\!<\!p\!\leq\!\infty\,,\,\alpha\!=\!1/q\!-\!1/p\,,$
we have for $0\!<\!t\!\leq\!1$
$$\sup_{a\in L^1(M\times N)}\,
{{\hbox{\nnorm{a}}_{p,q,t}}\over{\hbox{\Norm{a}}_{p,q,t}}} =
t^{-q/p-q}\,.$$
\par
\zz

\noi {\bf Proof.} Let $E,\,F$ with 
$t^{-1}\nu(F)^\alpha\!>\!\mu(E)^\alpha$
and
$$t\,\nu(F)^{-\alpha}\,
\Bigl(\int\!\!\!\int_{E\!\times\!F}\vert a(x,y)\vert^q\,dx\,dy\Bigr)
^{1/q}>1\,.
\leqno(3.2)$$

\noi In order to construct $\tilde E$ and $\tilde F\,,\,
t^{-1}\nu(\tilde F)^\alpha\!\leq\!\mu(\tilde E)^\alpha$ with
$$\mu(\tilde E)^{-\alpha}\,
\Bigl(\int\!\!\!\int_{\tilde E\!\times\!\tilde F}
\vert a(x,y)\vert^q\,dx\,dy\Bigr)
^{1/q}>t^{q/p-q}\,,                    
\leqno(3.3)$$

\noi we can suppose that 
$t^{-1}\nu(F)^\alpha\!>\!\mu(M)^\alpha$
and $E\!=\!M$ in (3.2).
If not, $\tilde F\!=\!F$ and any $\tilde E\!\supset\!E$ with 
$\mu(\tilde E)\!=\!t^{-1/\alpha}\,\nu(F)$ verify (3.3).

\hz

\noi By [BS], Theorem II, 2.7, p.\,51, 
we find for all $0\!<\!\lambda\!\leq\!1$
an $F_\lambda\!\subset\!N\,,~\nu(F_\lambda)\!=\!\lambda\,\nu(F)\,,$ such that
$$t\,\nu(F_\lambda)^{-\alpha}\,
\Bigl(\int\!\!\!\int_{M\!\times\!F_\lambda}
\vert a(x,y)\vert^q\,dx\,dy\Bigr)
^{1/q}>\lambda^{1/p}\,.$$

\noi If we put $\tilde E\!=\!M$ and $\tilde F\!=\!F_\lambda$ for
$$\lambda = {{\mu(M)}\over{\nu(F)}}\,t^{1/\alpha}\geq t^{1/\alpha}\,,$$

\noi then (3.2) is obvious.
\vfill
\eject
\noi
{\bf \S 4. Applications}
\dz

\noi
Let us denote by $B(X,Y)$ the space of bounded operators between two
Banach spaces $X,~Y$. Let $(M,\mu),~(N,\nu)$ be arbitrary
measure spaces. Let $B_0=B(L^\infty(d\nu), L^\infty(d\mu))$ and
$B_1=B(L^1(d\nu), L^1(d\mu))$. Clearly we may view $(B_0\,,B_1)$  as a
compatible couple in the sense of interpolation theory by identifying
an element of $B_0$ or $B_1$ with a linear operator from 
$L^\infty(d\nu)\cap L^1(d\nu)$ into $L^\infty(d\mu)+L^1(d\mu)$.
For simplicity we will assume (although this is inessential) that all
the spaces we consider are over the field of real scalars.

\hz
Let $L$, $\LL$ be real Banach lattices and let $u:~L\rightarrow \LL$ 
be a bounded operator. Then $u$ is called regular if there is a
positive operator $v:~L\rightarrow \LL$ such that
$$\vert u(x)\vert\leq v(\vert x\vert)~~~~~\hbox{for all }x~\hbox{in }L\,.$$

\noi
We define
$$\Vert u\Vert_r=\inf\{\Vert v\Vert\}$$

\noi
where the infimum runs over all such dominating operators $v$. We denote by
$B_r(L,\LL)$ the space of all regular operators equipped with this norm.
If $\LL$ is Dedekind complete in the sense of [MN] then $B_r(L,\LL)$
is a Banach lattice and we have simply
$$\Vert u\Vert_r=\Vert\,\vert u\vert\,\Vert_{\vrf{10}{0}B(L,\LL)}\,.$$

\noi
This applies in particular when $\LL$ is a Lorentz space, as in the 
situation we consider below. We should mention that any bounded 
operator between $L^\infty$--\,spaces (or  $L^1$--\,spaces) is automatically
regular, so that $B_0=B_r(L^\infty(d\nu),L^\infty(d\mu))$ and
$B_1=B_r(L^1(d\nu),L^1(d\mu))$.

\proclaim Theorem 4.1. For all $u$ in $B_0+B_1$, let
$$\hbox{\nnorm{u}}_t=\sup\,\Bigl\{\Big(\mu(E)\vee t^{-1}\nu(F)\Big)^{-1}
\Big\langle\vert u\vert(\char{F}),\char{E}\Big\rangle\Bigr\}$$
where the supremum runs over all measurable subsets 
$E\subset M,~F\subset N$ with positive and finite measure.
Then for all $t>0$
$${1\over 2}\,K_t(u\,;\,B_0\,,B_1)\leq
\hbox{\nnorm{u}}_t
\leq K_t(u\,;\,B_0\,,B_1)\,.$$
\par
\medskip
\noi{\bf Proof.} Assume first that $(M,\mu)$ and $(N,\nu)$ are purely
atomic measure spaces each with only finitely many atoms. Then $B_0$ and
$B_1$ can be identified (via their kernels) respectively with
$L^\infty(d\mu\,;\,L^1(d\nu))$ and $L^\infty(d\nu\,;\,L^1(d\mu))$.
Therefore, in that case the nontrivial part of Theorem 4.1 is but a
reformulation of Proposition 1.1.

\hz
\noi
The general case can be deduced from this using conditional expectations
and a simple weak compactness argument. We leave the details to the reader.

\zz
Let $(A_0,A_1)$ be any compatible couple of Banach spaces. We recall that
$(A_0,A_1)_{\theta,q}$ is defined for $0<\theta<1,~1\leq q\leq\infty$ as
the space of all $x$ in $A_0+A_1$ such that
$$\Vert x\Vert_{\theta,q} = 
\bigg(
\int_0^\infty\Big(t^{-\theta}K_t(x\,;\,A_0,A_1)\Big)^q\,{dt\over t}
\bigg)^{1/q}<\infty
$$

\noi
with the usual convention when $q=\infty$.
When equipped with the norm $\Vert\,.\,\Vert_{\theta,q}$, the space
$(A_0,A_1)_{\theta,q}$ is a Banach space. We refer to [BL] for more 
informations.

\zz
\proclaim Theorem 4.2. Let $0<\theta<1$ and $\theta=1/p\,.$ We have 
$$(B_0,B_1)_{\theta,\infty}= B_r\Big(L^{p,1}(d\nu),L^{p,\infty}(d\mu)\Big)$$
with equivalent norms.
\par
\medskip
\noi{\bf Proof.} 
Consider an operator $u:\,L^{p,1}(d\nu)\rightarrow L^{p,\infty}(d\mu)\,.$
Let us define for $\theta=1/p$
$$[\,u\,]_p=
\sup\,\bigg\{
\nu(F)^{-\theta}\mu(E)^{\theta-1}\,
\Big\vert\Big\langle u(\varphi\,\char{F})\,,
\psi\,\char{E}\Big\rangle\Big\vert\bigg\}
\leqno(4.1)$$
\noi
where the supremum runs over all $\varphi$ (resp. $\psi$) in the unit
ball of $L^\infty(d\nu)$ (resp. $L^\infty(d\mu)$) and over all 
measurable subsets $E\subset M,~F\subset N$ with finite positive measure.
By a well known property of the Lorentz spaces $L^{p,1}(d\nu)$ and
$L^{p^*,1}(d\mu)$ (cf.  [SW, Theorem 3.13]) there is a 
positive constant $C_p$
depending only on $p,~1<p<\infty$, such that we have
$$C_p^{-1}[\,u\,]_p\leq
\Vert u\Vert_{L^{p,1}(d\nu)\rightarrow L^{p,\infty}(d\mu)}\leq
C_p\,[\,u\,]_p\,,$$

\noi
for all $u:\,L^{p,1}(d\nu)\rightarrow L^{p,\infty}(d\mu)\,.$
If $u$ is positive, (4.1) can be simplified. In particular we have
$$[\,\vert u\vert\,]_p=
\sup_{E\subset M\atop F\subset N}
\,\bigg\{
\nu(F)^{-\theta}\mu(E)^{\theta-1}\,\Big\langle \vert u\vert(\char{F})\,,
\,\char{E}\Big\rangle\,\bigg\}
\leqno(4.2)$$

\noi
Assume $\nu(F)=t\,\mu(E)$ then 
$\nu(F)^{-\theta}\mu(E)^{\theta-1}=t^{-\theta}\mu(E)^{-1},$ hence
$$[\,\vert u\vert\,]_p=\sup_{t>0}\,t^{-\theta}\hbox{\nnorm{u}}_t\,.$$

\noi
Therefore by Theorem 4.1 $\Vert u\Vert_{(B_0,B_1)_{\theta,\infty}}$ is
equivalent to $[\,\vert u\vert\,]_p$ or equivalently to the norm
of $u$ in $B_r\Big(L^{p,1}(d\nu),L^{p,\infty}(d\mu)\Big)\,.$ This yields
immediately the announced result.

\hz

\noi {\bf Remark 4.1.} We refer to [P1] for the analogue of Theorem 4.2
for the complex interpolation method.

\goodbreak
\hz
We conclude this paper with an application to $H^p$--\,spaces in the framework
already considered in [P3]. Let $(A_0,A_1)$ be a compatible couple of 
Banach spaces and $S_0\subset A_0,$ $~S_1\subset A_1$ be closed subspaces.
Following the terminology in  [P2], we will say that $(S_0,S_1)$ is 
$K$--\,closed in $(A_0,A_1)$ if there is a constant $C$ such that for all
$x$ in $S_0+S_1$ and all $t>0$ we have
$$K_t(x\,;\,S_0,S_1)\leq C\,K_t(x\,;\,A_0,A_1)\,.$$

Let $(${\bf T}$, m)$ be the unit circle equipped with the normalized
Lebesgue measure. Let $B$ be a Banach space. We will denote for 
$1\leq p\leq\infty$ by $H^p(dm)$ (resp. $H^p(dm\,;\,B)$) 
the subspace
of $L^p(dm)$ (resp. $L^p(dm\,;\,B)$) of all the functions $f$
such that $\hat{f}(n)=0,~\forall\,n<0$. 

\hz
Let $(M,\mu)$ be any measure space. Consider the couple
$$X_0=L^1\Big(d\mu\,;\,L^\infty(dm)\Big),~~
X_1=L^1\Big(dm\,;\,L^\infty(d\mu)\Big)$$
\noi
and the subspaces
$$Y_0=L^1\Big(d\mu\,;\,H^\infty(dm)\Big),~~
Y_1=H^1\Big(dm\,;\,L^\infty(d\mu)\Big)\,.$$

\noi
It is proved in [P3, Lemma 2] that $(Y_0,Y_1)$ is $K$--\,closed in
$(X_0,X_1)$. By the simple duality principle emphasized in [P2]
(cf. Proposition 1.11 and Remark 1.12 in [P2]) this implies that
a similar property holds for the orthogonal subspaces $Y_0^\bot,~Y_1^\bot.$
More precisely, consider the subspaces
$$S_0=L^\infty\Big(d\mu\,;\,H^1(dm)\Big),~~
S_1=H^\infty\Big(dm\,;\,L^1(d\mu)\Big)$$
\noi
of the spaces
$$A_0=L^\infty\Big(d\mu\,;\,L^1(dm)\Big),~~
A_1=L^\infty\Big(dm\,;\,L^1(d\mu)\Big).$$

\noi
Then by this duality principle, $(S_0,S_1)$ is $K$--\,closed in $(A_0,A_1)$.
Therefore, our computation of the $K_t$--\,functional for the couple
$(A_0,A_1)$ (cf. Theorem 0.1 above) is applicable to the couple
$(S_0,S_1)$. Taking for simplicity $M={\bf N}$ equipped with the
counting measure, we obtain:

\hz
\proclaim Theorem 4.3. There is a numerical constant $C$ with the following
property. Let $(f_n)$ be a sequence in $H^1(dm)$ and let $t>0$.
Assume that for all subsets $E\subset{\bf N}$ and all measurable subsets
$F\subset{\bf T}$ we have
$$\sum_{n\in E}\,\int_F\vert f_n\vert\,dm\leq
\vert E\vert\vee t^{-1}m(F)\,.
\leqno(4.3)$$
\noi
Then there is a decomposition $f_n=g_n+h_n$ with 
$g_n\el H^1(dm),~h_n\el H^\infty(dm)$ such that
$$\sup_{n\in{\bf N}}\,\Vert g_n\Vert_{H^1(dm)}\leq C~~
\hbox{and}~~
\Big\|\sum_{n\in{\bf N}}\vert h_n\vert\,\Big\|_{L^\infty(dm)}
\leq C\,t^{-1}\,.$$

\vfill
\eject

\noi
{\bf \S 5. References}
\zz
\par
\frenchspacing
\begin{itemize}
\item[{[BF]}] R.\,C.\,Blei and J.\,J.\,F. Fournier. Mixed--norm conditions
and Lorentz norms. Proc. SLU--GTE Conference on Commutative Harmonic
Analysis {\flp 1987}. in:
Contemp. Math. {\bf 91} ({\flp 1989}). 57--78

\hz

\item[{[BL]}] J. Bergh and J. L\"ofstr\"om. 
     {\sl Interpolation Spaces: An Introduction}. Grundlehren der 
Mathematischen Wissenschaften {\bf 223}. Springer--Verlag. {\flp 1976}.

\hz

\item[{[BS]}] C. Bennett and R. Sharpley. 
     {\sl Interpolation of Operators}. Academic Press. {\flp 1988}.

\hz

\item[{[LT]}] J. Lindenstrauss and L. Tzafriri. 
     {\sl Classical Banach Spaces II}. Springer Verlag.  {\flp 1979}.

\hz

\item[{[MN]}] P. Meyer--Nieberg. 
     {\sl Banach Lattices}. Springer Verlag.  Univerisitext. {\flp 1991}.

\hz

\item[{[P1]}] G. Pisier. Complex interpolation and regular operators
     between Banach lattices. Archiv der Math. [\,to appear\,]

\hz

\item[{[P2]}] G. Pisier. Interpolation between $H^p$-spaces and 
     non-commutative generalizations I. 
     Pacific J. \,Math. {\bf 155} ({\flp 1992}). 341--368.

\hz

\item[{[P3]}] G. Pisier. Interpolation between $H^p$-spaces and 
     non-commutative generalizations II.
     Revista Mat. Iberoamericana  ({\flp 1993}) [\,to appear\,]

\hz

\item[{[SW]}] E. Stein and G. Weiss. ~
     {\sl Introduction to Fourier Analysis on Euclidean Spaces}. 
     Princeton University Press. {\flp 1971}.

\hz

\item[{[V]}] N. Varopoulos. On an inequality of von Neumann and an application 
     of the metric theory of tensor products to operators theory. 
     J. \,Funct.\,Anal. {\bf 16} ({\flp 1974}). 83--100.

\end{itemize}

\vz

\hbox{
\vtop{\hsize 7cm\obeylines
Gilles Pisier
\smallskip
Texas A.\,\&\,M. University
College Station
TX 77843~~--~~U.\,S.\,A.
\smallskip
and
\smallskip
Universit\'e Paris VI
Bo\^{\i}te 186
4, Place Jussieu
75252 Paris~~Cedex 05~~--~~France}
\hskip20pt
\vtop{\hsize 7cm\obeylines
Albrecht He\ss
\smallskip
Mathematische Fakult\"at
Friedrich--Schiller--Universit\"at
Uni--Hochhaus, 17.\,OG
07740 Jena~~--~~Germany}}

\vfill
\eject

\end{document}